\documentclass[a4paper,english,cleveref, autoref, thm-restate,
nolineno,
]{socg-lipics-v2021}

\pdfoutput=1 
\hideLIPIcs  

\usepackage{bm}
\usepackage{tikz-cd}
\usepackage{mathtools}
\usepackage{numprint}
\npthousandsep{,}
\npthousandthpartsep{}
\usepackage{booktabs}
\usepackage{multirow}
\usepackage{cellspace}
\usepackage{stackengine}
\usepackage{diagbox}
\usepackage{float}
\setstackEOL{\\}

\DeclareMathOperator{\SL}{SL}
\DeclareMathOperator{\PSL}{PSL}

\title{The complete 10-tetrahedra census of orientable cusped hyperbolic 3-manifolds} 

\author{Shana Yunsheng Li}{Dept. of Math., University of Illinois Urbana-Champaign, USA \and \url{https://shana-y-li.github.io} }{yl202@illinois.edu}{https://orcid.org/0009-0006-5320-106X}{}

\authorrunning{S. Y. Li} 

\Copyright{Shana Yunsheng Li} 



\category{} 

\supplementdetails[]{Dataset}{https://doi.org/10.7910/DVN/BFEXZA} 

\funding{This material is based upon work supported by the National Science Foundation under Grant No.\ DMS-2424139 while the author was in residence at the Simons Laufer Mathematical Sciences Institute in Berkeley, California, during the Spring 2026 semester. The author was also supported by the National Science Foundation under Grant No.\ DMS-2303572 during the same semester.}

\acknowledgements{The author wishes to thank Nathan Dunfield for his generous guidance during numerous discussions, and the Keeling cluster maintenance team at UIUC for providing the computational resources that made this project possible.}

\EventEditors{Hee-Kap Ahn, Michael Hoffmann, and Amir Nayyeri}
\EventNoEds{3}
\EventLongTitle{42nd International Symposium on Computational Geometry
(SoCG 2026)}
\EventShortTitle{SoCG 2026}
\EventAcronym{SoCG}
\EventYear{2026}
\EventDate{June 2--5, 2026}
\EventLocation{New Brunswick, NJ, USA}
\EventLogo{socg-logo.pdf}
\SeriesVolume{367}
\ArticleNo{XX}     

\begin{document}

\maketitle

\begin{abstract}
We extend the complete census of orientable cusped hyperbolic $3$-manifolds to $10$ tetrahedra, giving the next $\numprint{150730}$ manifolds and their $\numprint{496638}$ minimal ideal triangulations. As applications, we find the precisely $\numprint{439898}$ exceptional Dehn fillings on them, revealing the next $\numprint{1849}$ simplest hyperbolic knot exteriors in $S^3$. We also give the simplest example of an orientable cusped hyperbolic $3$-manifold containing a closed totally geodesic surface. 
\end{abstract}

\section{Introduction}
\label{sec:intro}

The enumeration of examples has been of fundamental importance in the study of low-dimensional topology ever since Tait gave the first knot table in 1884 \cite{Tait:knot}. The censuses of knots have been extensively explored over the years, all the way to approximately the first 2 billion knots \cite{HTW:knot,Burton:knot,Thistlethwaite:knot}. However, the censuses of orientable cusped hyperbolic $3$-manifolds, initiated by Hildebrand and Weeks \cite{HW:census}, have seen no progress after the $\numprint{44250}$ manifolds whose minimal ideal triangulations consist of $9$ tetrahedra were given by Burton in 2014 \cite{BurtonBenjaminA2017Tchc}. In this paper, we extend the census by giving the $\numprint{150730}$ manifolds whose minimal ideal triangulations consist of $10$ tetrahedra. \Cref{tab:censuses} summarizes the timeline of the development of censuses of orientable cusped hyperbolic $3$-manifolds.

\begin{table}[ht]
  \centering
  \begin{tabular}{cccc}
            \toprule
            Tetrahedra & Name of manifolds & Year & Contributor(s) \\
             \midrule
               2 - 5 & m003 $\sim$ m412 & 1989 & Hildebrand \& Weeks \cite{HW:census}\\
            \midrule
            6 & s000 $\sim$ s961 & \multirow{2}{*}{1999} & 
            \multirow{2}{*}{\Centerstack{Callahan, \\ Hildebrand \& Weeks \cite{CHW:census}} } \\
            7 & v0000 $\sim$ v3551  & & \\
            \midrule
            8 & t00000 $\sim$ t12845 & 2010 & Thistlethwaite \cite{Thistlethwaite:8-tet} \\
            \midrule
            9 & o9\_00000 $\sim$ o9\_44249 & 2014 & Burton \cite{BurtonBenjaminA2017Tchc}\\
            \midrule
            10 & o10\_000000 $\sim$ o10\_150729 & 2025 & L.\\
            \bottomrule
  \end{tabular}
  \caption{Timeline of censuses of orientable cusped hyperbolic $3$-manifolds}
  \label{tab:censuses}
  \end{table}


By hyperbolic $3$-manifolds, we mean $3$-manifolds with Riemannian metrics of constant curvature $-1$ that are complete as metric spaces; in particular their geodesics can be extended indefinitely in length. A hyperbolic $3$-manifold is \textit{cusped} if it has a topological end that is of finite volume. An orientable cusped finite-volume hyperbolic $3$-manifold is homeomorphic to the interior of a compact $3$-manifold whose boundary is a non-empty union of tori. Our main result is the following:

\begin{theorem}
  \label{th:10census}
  There are precisely $\numprint{150730}$ orientable cusped hyperbolic $3$-manifolds whose minimal ideal triangulations consist of $10$ tetrahedra. Moreover, there are precisely a total of $\numprint{496638}$ minimal ideal triangulations of these manifolds.
\end{theorem}
\Cref{tab:cusps} in \cref{subappen:cusps} gives the counts of manifolds in \cref{th:10census} by cusp number.

We highlight here that in order to obtain \cref{th:10census}, we applied a new core technique, the verified canonical triangulations. We anticipate that a mixture of this technique with traditional ones would suffice for further extension of the census, provided successful generation of the candidates; see \cref{subsec:overview,subsec:obstruct} for more detailed discussions.

Throughout this paper, we assume triangulations consist of finitely many tetrahedra.
If a cusped hyperbolic $3$-manifold admits an ideal triangulation, 
each tetrahedron has finite volume \cite[\S 2.5]{Dunfield:volume}, hence the volume of the manifold is finite. 
Conversely, any cusped finite-volume hyperbolic $3$-manifold has an ideal triangulation given by a refinement of the Epstein--Penner cellulation \cite[\S 5]{Martelli:intro}.
Therefore, for cusped hyperbolic $3$-manifolds, being finite-volume is equivalent to admitting ideal triangulations, hence we drop the term finite-volume when talking about cusped hyperbolic $3$-manifolds admitting ideal triangulations.

In general, it is unknown whether all cusped finite-volume hyperbolic $3$-manifolds admit geometric triangulations in the sense of \cref{subsec:triangulation} below
\cite{LST:virtual}. However, when certifying hyperbolicity, \texttt{SnapPy} returns true if and only if it finds a geometric triangulation.
Therefore an immediate corollary of this paper is that 
\begin{corollary}
  \label{cor:geo}
 All cusped hyperbolic $3$-manifolds ideally triangulable by at most $10$ tetrahedra admit geometric triangulations. 
\end{corollary}
Note that nonorientable manifolds are included here since they are included in \cref{th:candidates}. For orientable manifolds, the above result can be strengthened to \cref{pr:min-geo}.

We briefly review the needed terminologies in \cref{sec:prelim}, then give detailed discussions in \cref{sec:creation} for the creation of the $10$-tetrahedra census, after which we give some applications of the extended census in \cref{sec:application}. All corresponding code and data are available at \cite{10-tets:dataverse}.

\section{Preliminaries}
\label{sec:prelim}

\subsection{Triangulations}

\label{subsec:triangulation}

A \textit{(generalised) triangulation} $\mathcal{T}$ is a disjoint union of finitely many tetrahedra, with some or all faces identified in pairs via simplicial maps, and we avoid identifications where an edge is identified with itself in reverse. Note that $\mathcal{T}$ is not necessarily a simplicial complex, as identifications of two faces of the same tetrahedron are allowed. 

Quotienting by the identifications, $\mathcal{T}$ can be seen as a topological space, along with a cellular structure induced by the tetrahedra. A \textit{vertex} in $\mathcal{T}$ is a $0$-cell in this cellular structure, and an \textit{edge} is a $1$-cell. The \textit{link} of a vertex in $\mathcal{T}$ is the frontier of a small regular neighborhood of that vertex. A vertex is \textit{internal} if its link is a sphere, 
and is \textit{ideal} if its link is some other closed surface. For an edge $e$ in $\mathcal{T}$, its \textit{degree} is the cardinality of the set
\begin{equation*}
  \mathcal{S}_e \coloneqq \bigcup_{t \in \{3\text{-cells of }\mathcal{T}\}} \{\text{(nonempty) connected components of }\mathcal{N}_e\cap t\},
\end{equation*}
where $\mathcal{N}_e$ is a small regular neighborhood of the midpoint of $e$.

\begin{definition}
  \label{def:ideal-triangulation}
  Let $M$ be the interior of a compact $3$-manifold with boundary. A triangulation $\mathcal{T}$ is an \textit{ideal triangulation of $M$} if it is homeomorphic to $M$ after removing its ideal vertices, and is \textit{minimal} if no ideal triangulation of $M$ with fewer tetrahedra exists.
\end{definition}
It follows that every vertex of an ideal triangulation is either internal or ideal. Note that an ideal triangulation itself is not a manifold and has at least one ideal vertex. 

If $M$ is a cusped hyperbolic $3$-manifold, an ideal triangulation $\mathcal{T}$ of $M$ is \textit{geometric} if all vertices are ideal, all edges are geodesics, all $2$-cells are parts of geodesic planes, and no flat or negatively oriented tetrahedron (i.e. with zero or negative volume) exists. 

The following results were proven in Section 3 of \cite{BurtonBenjaminA2017Tchc}:

\begin{lemma}[Burton]
  \label{le:no-internal}
  Let $\mathcal{T}$ be a minimal ideal triangulation of a cusped hyperbolic $3$-manifold. Then $\mathcal{T}$ contains no internal vertices.
\end{lemma}

\begin{lemma}[Burton]
  \label{le:degree}
  Let $\mathcal{T}$ be a minimal ideal triangulation of a cusped hyperbolic $3$-manifold. Then $\mathcal{T}$ has no edges of degree $1$ or $2$, and $\mathcal{T}$ has no edges of degree $3$ that are contained in three distinct tetrahedra.
\end{lemma}

\subsection{Normal surfaces}

In a tetrahedron, an arc $\alpha$ properly embedded in one of its faces is \textit{normal} if the endpoints of $\alpha$ are in the interior of distinct edges. A properly embedded disk is a \textit{triangle} if its boundary consists of three normal arcs in distinct faces, and is a \textit{square} if its boundary consists of four normal arcs in distinct faces. 

Let $\mathcal{T}$ be a triangulation. 
A properly embedded compact surface in $\mathcal{T}$ is \textit{normal} if its intersection with every tetrahedron in $\mathcal{T}$ (with the identifications undone) is a properly embedded surface whose connected components consist of only triangles and squares.

An isotopy of $\mathcal{T}$ is normal if it preserves the cellular structure of $\mathcal{T}$.
Up to normal isotopy, there are $4$ triangles and $3$ squares in a tetrahedron. Multiple copies of triangles or squares in the same normal isotopy class may appear in the intersection of a normal surface and a tetrahedron. It follows that a normal surface in a triangulation of $n$ tetrahedra is uniquely determined by a vector in $\mathbb{N}^{7n}$ up to normal isotopy, where $\mathbb{N}\coloneqq \{k\in\mathbb{Z}\mid k\geq 0\}$. 

The subset $\mathcal{N} \coloneqq \{\bm v\in \mathbb{N}^{7n}\mid \bm v\text{ represents a normal surface}\}$ consists of integral vectors in a certain non-convex subset of a polyhedral cone in $\mathbb{R}^{7n}$ \cite{JO:haken}. The \textit{standard vertex normal surfaces} correspond to vertices of a cross-section of this polyhedral cone, and can be viewed as forming a finite basis for $\mathcal{N}$.

The following result was proven in Section 3 of \cite{BurtonBenjaminA2017Tchc}:

\begin{lemma}[Burton]
  \label{le:normal}
  Let $\mathcal{T}$ be a minimal ideal triangulation of a cusped hyperbolic $3$-manifold. Then no standard vertex normal surface in $\mathcal{T}$ 
  \begin{itemize}
    \item has positive Euler characteristic; or
    \item has zero Euler characteristic, is not a vertex link, and in addition
    \begin{itemize}
      \item is one-sided; or
      \item is two-sided, and if we cut $\mathcal{T}$ open along it then no component of the resulting space is a solid torus or solid Klein bottle.
    \end{itemize}
  \end{itemize}
\end{lemma}

\section{Creating the 10-tetrahedra census}
\label{sec:creation}

\subsection{Overview}
\label{subsec:overview}

To build a census (of knots, manifolds, etc), the procedure generally follows a $3$-step outline:
\begin{romanenumerate}
    \item Generate all candidates for the census.
    \item Decide the eligibility of each candidate.
    \item Group eligible candidates into isomorphism classes and certify that the isomorphism classes are all distinct.
\end{romanenumerate}
As long as the candidate generation is performed exhaustively and the eligibility for each candidate is decided correctly, the census obtained is guaranteed to contain all wanted objects. With step (iii), the resulting census will contain no duplicates. We say that a census is \textit{complete} if it contains all wanted objects without any duplicated ones.

For censuses of cusped hyperbolic $3$-manifolds, the candidates will be ideal triangulations consisting of a given number of tetrahedra, and the eligibility conditions will be minimality (of triangulations) and hyperbolicity (of manifolds the triangulations represent). For certifying distinctness, algebraic invariants were widely applied for both the census of hyperbolic $3$-manifolds \cite{BurtonBenjaminA2017Tchc} and the census of knots \cite{HTW:knot,Burton:knot,Thistlethwaite:knot}. However, it is unknown whether algebraic methods will indeed separate all given isometry classes within a given scope, consequently they cannot be used to certify if two isometry classes are indeed the same, and other methods will be required for that purpose. Here we propose an alternative approach for cusped hyperbolic $3$-manifolds, using a geometric invariant, the \textit{canonical triangulation}, built upon the Epstein--Penner canonical cellulation \cite{EpsteinPenner:cell} by Weeks \cite{Weeks:convex}, whose rigorous computation was implemented by Goerner \cite{SnapPy:Verified}. An introduction to these can be found in \cite{FGGTV:tetrahedral}, and we review them briefly in \cref{subsubsec:canonical}. The canonical triangulation, when computed rigorously, is a complete invariant of cusped hyperbolic $3$-manifolds, and we applied it to deduplicate a majority of triangulations in the $10$-tetrahedra census in \cref{subsubsec:exact}.

Our approaches for steps (i) and (ii) largely overlap with those in \cite{BurtonBenjaminA2017Tchc}. For completeness, we briefly describe them in \cref{subsec:test}, and kindly refer the reader to \cite{BurtonBenjaminA2017Tchc} for more details. For step (iii), detailed discussions will be given in \cref{subsec:dedup}.

\subsection{Candidate generation and eligibility tests}
\label{subsec:test}

Let us start by looking at \cref{tab:steps}. The steps proceed from top to bottom. Each row consists of the name of the step, the number of remaining candidates whose hyperbolicity is unknown, the numbers of candidate discarded within that step due to ineligibility, and the number of candidates kept to be dealt with later in \cref{subsubsec:exact}.

\begin{table}[ht]
  \centering
  \begin{tabular}{crrr}
            \toprule
            Step & \multicolumn{1}{c}{Candidates} & \multicolumn{1}{c}{Discarded} & \multicolumn{1}{c}{Kept} \\
             \midrule
             \midrule
               Candidate generation & \numprint{8373308} & 0 & 0 \\
            \midrule
            Greedy nonminimality & \numprint{1946782} & \numprint{6426526} & 0 \\
            \midrule
            Certify hyperbolicity 
            ($r = 1$)  & \numprint{1072874} & 0 & \numprint{873908} \\
            \midrule
            Exhaustive nonminimality ($h = 2$) & \numprint{698650} &\numprint{374224} & 0 \\
            \midrule
            Special surfaces & \numprint{33807} & \numprint{664843} & 0 \\
            \midrule
           Certify hyperbolicity  ($r = 60$)  &\numprint{3269} & 0 &\numprint{30538}\\ 
            \midrule
            Exhaustive nonminimality ($h = 5$) & 12 & \numprint{3257} & 0\\
            \midrule
            Exhaustive nonminimality 
            ($h = 6$) & 6 & 6 & 0 \\
            \midrule
            Certify hyperbolicity 
            ($r = 1000$)  & 0 & 0 & 6\\
            \midrule
            \midrule
            Total  & \numprint{8373308} & \numprint{7468856} & \numprint{904452}\\
            \bottomrule
    \end{tabular}
    \caption{The steps of filtering candidates for the $10$-tetrahedra census}
    \label{tab:steps}
\end{table}

For candidate generation, we used \texttt{Regina}'s utility command \texttt{tricensus} to generate all ideal triangulations consisting of $10$ tetrahedra except those that are obviously nonminimal or non-hyperbolic. More precisely, we obtained:

\begin{theorem}
    \label{th:candidates}
    Up to relabelling, there are $\numprint{8373308}$ ideal triangulations with $10$ tetrahedra, in which every vertex is ideal, there are no edges of degree $1$ or $2$, and there are no edges of degree $3$ that belong to three distinct tetrahedra.  
\end{theorem}

Due to \cref{le:no-internal,le:degree}, the triangulations described in \cref{th:candidates} contain all minimal ideal triangulations of cusped hyperbolic $3$-manifolds whose minimal ideal triangulations consist of $10$ tetrahedra. Note that triangulations representing nonorientable manifolds are also included and carried through the tests, though we will discard them in \cref{subsec:dedup}. 

With the candidates generated, we proceed to determine their eligibilities. For an ideal triangulation to remain in the census, it must be both minimal and representing a cusped hyperbolic $3$-manifold. Hence an ideal triangulation is discarded immediately if it is certified to be nonminimal or representing a non-hyperbolic $3$-manifold. If an ideal triangulation is certified to represent a cusped hyperbolic $3$-manifold, we temporarily keep it, and decide whether it is minimal by checking if it represents a manifold that has already appeared in censuses with fewer tetrahedra, which we discuss in \cref{subsec:dedup}. 

Four tests were used in \cref{tab:steps}, which we describe below:
\begin{description}
  \item[Greedy nonminimality:]  we use \texttt{Regina}'s highly efficient simplification algorithm \cite{BurtonRegina} to try converting each candidate into an ideal triangulation with fewer tetrahedra, while representing the same $3$-manifold. If such a triangulation with fewer tetrahedra is found, then the candidate is nonminimal. 
  \item[Exhaustive nonminimality:] given an integer parameter $h>0$, for each candidate with $n$ tetrahedra ($n=10$ in our settings), enumerate all triangulations that can be reached via $2$-$3$ and $3$-$2$ Pachner moves without ever exceeding $n+h$ tetrahedra. If a triangulation with fewer than $n$ tetrahedra is found during the process, the candidate is nonminimal. 
  \item[Special surfaces:] for each candidate $\mathcal{T}$, enumerate all standard vertex normal surfaces in it, and check if any of them
  \begin{itemize}
    \item has positive Euler characteristic; or
    \item has zero Euler characteristic, is not a vertex link, and in addition
    \begin{itemize}
      \item is one-sided; or
      \item is two-sided, and if we cut $\mathcal{T}$ open along it then no component of the resulting space is a solid torus or solid Klein bottle.
    \end{itemize}
  \end{itemize}
  If any such surface is found, then by \cref{le:normal}, $\mathcal{T}$ is either nonminimal or does not represent a cusped hyperbolic $3$-manifold.
  \item[Certify hyperbolicity:] given an integer parameter $r>0$, for each candidate $\mathcal{T}$, we randomize it by doing random Pachner moves, and ask \texttt{SnapPy} to try finding a complete hyperbolic structure on $\mathcal{T}$, which uses an algorithm similar to HIKMOT's \cite{HIKMOT}, where affirmative answers are rigorous \cite{SnapPy:Verified}. If the answer is not affirmative, we randomize $\mathcal{T}$ and try again, until an affirmative answer is achieved or the number of attempts reaches $r$. 
\end{description}
We summarize the tests along with their core functions in \cref{tab:tests} in \cref{subappen:summary}.

As is presented in \cref{tab:steps}, for each candidate, an affirmative answer was given by one of the above four tests. A candidate was kept if and only if the affirmative answer was given by the certify hyperbolicity test, and was discarded otherwise.
As such, we obtained $\numprint{904452}$ triangulations certified to represent cusped hyperbolic $3$-manifolds at the bottom of \cref{tab:steps}, with all others discarded. 


\subsection{Grouping and Separating}
\label{subsec:dedup}



In this section, we discuss the procedure for step (iii) described in \cref{subsec:overview}. A central tool for it, in both \cite{BurtonBenjaminA2017Tchc} and our method, is the \textit{canonical triangulation} of cusped finite-volume hyperbolic $3$-manifolds. 


\subsubsection{Canonical triangulation}
\label{subsubsec:canonical}

We review the definition and computation of canonical triangulations very briefly, and refer the reader to \cite{FGGTV:tetrahedral} for more details. 

Consider the $(3+1)$-Minkowski space $\mathbb{E}^{3,1} \coloneqq \mathbb{R}^4$ whose inner product $\langle\cdot,\cdot\rangle$ is defined by $\langle \bm x, \bm y\rangle \coloneqq \sum_{i=0}^2x_iy_i  - x_3y_3$. The hyperboloid $S^+\coloneqq \{\bm x\in \mathbb{E}^{3,1} \mid x_3>0,\quad \langle \bm x,\bm x\rangle = -1 \}$ is isometric to $\mathbb{H}^3$. Given a cusped hyperbolic $3$-manifold $M$, for each cusp of $M$ we choose a horotorus cusp neighborhood of the same volume; each neighborhood is bounded by a torus, which lifts to horospheres in $\mathbb{H}^3\cong S^+ $. For each horosphere $\mathcal{H}\subset S^+$, its normal vector $\bm v$ in $\mathbb{E}^{3,1}$ satisfies $\langle \bm v ,\bm v \rangle = 0$, and we normalize it so that $\langle \bm v,\bm w\rangle = -1$ for any $\bm w\in \mathcal{H}$. Therefore we obtain a family $\{ \bm v\}_M $ of such vectors associated to a cusped hyperbolic $3$-manifold $M$.

\begin{definition}
  The Epstein--Penner cellulation of $M$ is given by the radial projection of the polygonal faces of the convex hull of $\{\bm v\}_M $ onto $S^+$.
\end{definition}

The Epstein--Penner cellulation is invariant under the choice of the volume of cusp neighborhoods, hence it is an invariant of cusped hyperbolic $3$-manifolds \cite{EpsteinPenner:cell}. It is generically a triangulation, but not always so. Given a polygonal cellulation of a $3$-manifold, one can subdivide it into the triangulation induced by the following refinement of its $1$-skeleton:
\begin{itemize}
  \item add a vertex in the interior of each $2$-cell and $3$-cell;
  \item for each $2$-cell or $3$-cell $e$, let $v$ be the vertex added to its interior:
  \begin{itemize}
    \item for each vertex $u$ originally in $e$, add an edge connecting $v$ and $u$;
    \item if $e$ is a $3$-cell, for each $2$-cell $e_2$ on the boundary of $e$, add an edge connecting $v$ and the vertex added to the interior of $e_2$.
  \end{itemize} 
\end{itemize}

\begin{definition}
  The canonical triangulation of $M$ is the Epstein--Penner cellulation of $M$ if the latter is a triangulation, and the subdivision of the Epstein--Penner cellulation described above otherwise.
\end{definition}

Weeks established a heuristic procedure which, given a geometric triangulation, 
computes the Epstein--Penner cellulation of the corresponding cusped hyperbolic $3$-manifold via local modifications, doing $2$-$3$ and $3$-$2$ Pachner moves wherever tetrahedra meet concavely at faces or edges \cite{Weeks:convex}. The canonical triangulation can then be easily computed from the Epstein--Penner cellulation. The methods for both were implemented in \texttt{SnapPea} in the 1990s, and later inherited by its successor \texttt{SnapPy} \cite{SnapPy}. In their original implementations, the hyperbolic structures used were only recorded numerically, hence accuracy issues may come into play, which we discuss in the following two sections. 

\subsubsection{Burton's method: algebraic certification}
\label{subsubsec:burton}

Let us briefly review the method Burton used in \cite{BurtonBenjaminA2017Tchc}, which can be summarized as:
\begin{enumerate}
  \item For each ideal triangulation, compute its canonical triangulation numerically using \texttt{SnapPea}, and partition them into classes by the canonical triangulations computed.
  \item For triangulations in each class, perform an exhaustive search through $2$-$3$ and $3$-$2$ Pachner moves, attempting to identify classes that should have been merged in step 1, but were not due to numerical accuracy issues, thus refining the partitioning.
  \item For each class, compute the list $\{\text{orientability}, H_1 = \pi_1^{\text{ab}}, S_2^{\text{ab}}(\pi_1),\dots, S_{k}^{\text{ab}}(\pi_1)\}$ of algebraic invariants to certify that all groups are distinct, increasing $k$ if necessary, where $S_i^{\text{ab}}(G)$ 
     is the multi-set of abelianisations of index $i$ subgroups of $G$, one subgroup for each conjugacy class. For the $9$-tetrahedra census, $k = 11$ was required and sufficed.
\end{enumerate}
Note that to confirm the completeness of previous censuses, Burton dealt with all candidates with no more than $9$ tetrahedra at once, hence triangulations that were grouped with triangulations with fewer tetrahedra were discarded immediately due to nonminimality.

There are several aspects where numerical computation could cause issues in step 1:
\begin{itemize}
  \item A quantity \textit{tilt} is used to characterize whether tetrahedra meet concavely, parallelly or convexly, each corresponding to the tilt being positive, zero, or negative respectively \cite{Weeks:convex,SW:tilt}. When the absolute value of tilt is small, the numerical computation may falsely recognize the sign of tilt, resulting in an incorrect canonical triangulation.
  \item Flat tetrahedra being falsely recognized as non-flat, which would then be incorrectly dealt with when subdividing Epstein--Penner triangulations into canonical triangulations. 
\end{itemize}

When the above issues arise in step 1, \texttt{SnapPea}/\texttt{SnapPy} fails to find the correct canonical triangulation and arrives at an incorrect one, while still representing the same manifold as the input triangulation. Hence steps 2 and 3 were necessary to spot and resolve these miscomputed cases. However, steps 2 and 3 can be very expensive to perform, and require back-and-forth adjustments of parameters. Furthermore, although successful for censuses with at most $9$ tetrahedra, it is unknown whether step 3 will separate all cases within a computable scope, or even if it will separate all of them in theory \cite{Gardam:profinite}.

\begin{remark}
  We remark here that the subdividing algorithm Weeks implemented in \texttt{SnapPea}/\texttt{SnapPy} does not work for refinements of Epstein--Penner cellulations containing flat tetrahedra. This led to the miscomputed example x101 
  found by Burton \cite{Burton:duplicate} and mentioned in \cite[Remark 3.6]{FGGTV:tetrahedral}. 
  The cause of the miscomputation was not located at the time of \cite{FGGTV:tetrahedral}. During the course of this project, the author found another such example, o9\_43806
  and located the said cause. This bug has been fixed in \texttt{SnapPy 3.3} by forcing it to find a refinement that does not contain flat tetrahedra before proceeding to subdividing it. After the fix, the numerical computation gives correct results for all minimal triangulations in the orientable cusped censuses of at most $10$ tetrahedra, which we confirmed using the verified computation explained in \cref{subsubsec:exact} (as long as successfully computed).
\end{remark}

\subsubsection{Our method: verified computation}
\label{subsubsec:exact}

Instead of making up for the mistakes caused by numerical computation of canonical triangulations using other methods, we enhance the computation itself to avoid the mistakes directly. With the accuracy issues described in the previous section, this requires us to:
\begin{alphaenumerate}
  \item given two supposedly different quantities, rigorously tell which one is larger;
  \item given two supposedly equal quantities, prove that they are equal.
\end{alphaenumerate} 

For (a), we use interval arithmetics, which can be thought of as computation with controlled errors. If two quantities are different, with sufficiently small error the lower bound of one will be strictly larger than the upper bound of another, revealing the larger quantity. Similar approach was also used in HIKMOT's rigorous certification of hyperbolicity \cite{HIKMOT}.

For (b), we use the LLL algorithm \cite{LLL:algorithm} to symbolically solve the trace field of hyperbolic structure, representing it as a number field. We then express the quantities as elements in the number field, hence they are equal if and only if their difference is the zero element.

The above enhanced computation was implemented in \texttt{SnapPy} by M. Goerner \cite{SnapPy:Verified}, and is called verified computation. We refer the reader to \cite{DHL:asymmetric,FGGTV:tetrahedral} for more detailed discussions. 

Canonical triangulations obtained with verified computation are thus provably correct. Let us call a supposedly canonical triangulation obtained by \texttt{SnapPy} a \textit{verified canonical triangulation} if it was obtained with verified computation, and an \textit{unverified canonical triangulation} if the computation was performed purely numerically. 

Let $\{\mathcal{T}\}$ be the set of ideal triangulations representing cusped hyperbolic $3$-manifolds and $\{\mathcal{M}\}$ be the set of cusped finite-volume hyperbolic $3$-manifolds. 
Let $\iota\colon \{\mathcal{T}\}\to \{\mathcal{M}\}$ be the function that maps an ideal triangulation to the manifold it represents. The unverified canonical triangulation can thus be seen as a multi-valued function $c_u\colon \{\mathcal{M}\} \to \{\mathcal{T}\}$ which is a right inverse of $\iota$; in other words, the following diagram commutes:
\begin{equation}
\begin{tikzcd}
  & \{\mathcal{T}\} \arrow[r, "\iota"] & \{\mathcal{M}\}\\
  &  &  \{\mathcal{M}\}\arrow[ul, bend left = 10, "c_u"]\arrow[u,swap, "1_{\{\mathcal{M}\}}"]
\end{tikzcd}
\end{equation}

In particular, $c_u$ is injective. The verified canonical triangulation is then a branch of $c_u$, a single-valued function $c_v\colon \{\mathcal{M}\} \to \{\mathcal{T}\}$. Since $c_u$ is injective, $c_v$ is also injective, hence restricts to a bijection between $\{\mathcal{M}\}$ and its image. Therefore the verified canonical triangulation is a complete invariant of cusped finite-volume hyperbolic $3$-manifolds.

Although the definition of canonical triangulation does not require the manifold to be orientable, the current implementation of verified computation for canonical triangulations in \texttt{SnapPy} does not work for nonorientable ones with multiple cusps.
Hence we discard the nonorientable candidates and deal with the remaining $\numprint{608918}$ orientable ones.

We successfully computed the verified canonical triangulations for all except $16$ of the $\numprint{608918}$ orientable candidates. The $16$ cases failed due to their trace fields being too complicated; we remark that this very nature intrinsically separates them from the others, though further investigation is still required to prove it. By comparing the results with the verified canonical triangulations of all manifolds in censuses of at most $9$ tetrahedra, we obtained $\numprint{496622}$ minimal triangulations of $10$ tetrahedra, which, grouped by their verified canonical triangulations, give $\numprint{150724}$ distinct manifolds. For the $16$ candidates whose verified canonical triangulations were not computed, we computed the first $60$ digits of their hyperbolic volumes using interval arithmetics \cite{SnapPy:Verified} and their first homology groups, which separated them from all manifolds in censuses up to $9$ tetrahedra plus the $\numprint{150724}$ ones, and grouped them into $6$ classes; within each class, the triangulations were confirmed to all represent the same manifold as their unverified canonical triangulations are the same.

Therefore, we conclude that there are exactly $\numprint{150730}$ isometry classes in the $10$-tetrahedra census, containing a total of $\numprint{496638}$ minimal triangulations, obtaining \cref{th:10census}.

\subsection{Finalizing the census}
\label{subsec:finalize}

For each of the $\numprint{150730}$ isometry classes obtained in \cref{subsec:dedup}, we select a representative triangulation by the following procedure:
\begin{bracketenumerate}
\item If the verified canonical triangulation is in the isometry class, let it be the representative.
\item Otherwise, consider the subset of triangulations that are geometric, which turns out to be never empty. 
If there is only one triangulation in the class, let it be the representative.
\item If there are multiple triangulations in that subset, let the triangulation with the most self-isomorphisms, i.e. cellular homeomorphisms onto itself, be the representative.
\item If there are multiple triangulations with the most self-isomorphisms, choose the one with the ``fattest'' tetrahedra, i.e. the one with the largest minimum imaginary part of tetrahedron shape parameters.
\item If there are still multiple options left, choose the one with the smallest triangulation signature (as defined in \cite{Burton:signature}) under the lexicographical order. 
\end{bracketenumerate}

With the representative triangulations chosen as above, we set their peripheral curves to be the shortest. We then ordered them by volume, named them as displayed in \cref{tab:censuses}, saved and published the data in \cite{10-tets:dataverse}. The data now comes automatically with \texttt{SnapPy 3.3}. 

Note that verified canonical triangulations, by construction, are always geometric 
\cite{FGGTV:tetrahedral}. Hence due to the geometric requirement in step (2), we obtained
\begin{proposition}
  \label{pr:min-geo}
  Every cusped hyperbolic $3$-manifold in the $10$-tetrahedra census admits a minimal triangulation that is geometric. 
\end{proposition}
We remark that the above proposition also holds for all (orientable) censuses with fewer tetrahedra, which we checked using Burton's data from \cite{BurtonBenjaminA2017Tchc}.

\subsection{Obstructions on further extending the census}
\label{subsec:obstruct}

As is discussed in \cref{subsec:overview}, there are essentially three steps for creating a complete census. For the census of (orientable) cusped hyperbolic $3$-manifolds, step (ii) should not be a problem, as there are several other tests proposed in \cite{BurtonBenjaminA2017Tchc}, but were not used when creating the $10$-tetrahedra census. For step (iii), the verified computation has a chance of failing only when the LLL algorithm is invoked, which happens only when a tilt that is extremely close to zero is encountered during computation. This essentially implies that the Epstein--Penner cellulation is not a triangulation, which is the case for only $986$ triangulations out of the $608918$ orientable candidates mentioned in \cref{subsubsec:exact}, a frequency less than $0.2\%$. Therefore, it is promising that a combination of verified canonical triangulations with algebraic methods would suffice in the foreseeable future, where the former separates most manifolds, and the latter deals with the rest.

The main obstruction thus lives in step (i), where an exponential growth rate of the computational power is needed, with a much larger base compared to tabulating knots \cite{RFS:alt}. If one were to include nonorientable candidates during generation (as we did in \cref{th:candidates}), even though the algorithm established in \cite{BurtonBenjaminA2017Tchc} is already extremely efficient, it still takes approximately $2$ years to generate all candidates with $10$ tetrahedra in a single thread of $2.6$ GHz, about $50$ times longer than to generate all candidates with 9 tetrahedra; we accelerated the process by parallelizing it in $300$ threads. Thus, assuming the same growth rate, it will take approximately $100$ years to generate all candidates with $11$ tetrahedra with a single thread of $2.6$ GHz, and hence $4$ months in $300$ threads. 

The time consumption and growth rate can both be lowered by excluding nonorientable candidates during generation, as half of identifications of each face will be excluded from the beginning. It takes $2$ months to enumerate all orientable candidates with $10$ tetrahedra in a single thread, about $25$ times longer than for $9$ tetrahedra. We successfully generated all orientable candidates with $11$ tetrahedra, which took about $6$ years in a single thread, giving \cref{th:candidates-11}. We have now completed the $11$-tetrahedra orientable census, with a preliminary \texttt{SnapPy} module available at \cite{snappy_11_tets}, and full details to appear in a future publication.

\begin{theorem}
    \label{th:candidates-11}
    Up to relabelling, there are $\numprint{27794289}$ ideal triangulations of orientable $3$-manifolds with $11$ tetrahedra, in which every vertex is ideal, there are no edges of degree $1$ or $2$, and there are no edges of degree $3$ that belong to three distinct tetrahedra.  
\end{theorem}

\section{Applications}
\label{sec:application}

A substantial amount of follow-up work can be done with the extended census. We describe several of them in this section. 

\subsection{Extending the census of exceptional Dehn fillings}
\label{subsec:exceptional}

Let $\overline{M}$ be a compact orientable $3$-manifold whose interior $M$ is a $1$-cusped hyperbolic $3$-manifold; in particular $\partial \overline{M}$ is a torus, and $M$ has finite volume. Dehn fillings on $M$ are characterized by slopes on $\partial \overline{M}$, which are unoriented isotopy classes of simple closed curves on the torus. For a slope $\alpha$, we denote by $M(\alpha)$ the closed manifold obtained by Dehn filling on $M$ so that $\alpha$ bounds a disk in the attached solid torus $D^2\times S^1$. 
\begin{definition}
  A slope $\alpha$ is exceptional if $M(\alpha)$ is not hyperbolic.
\end{definition}
The census of exceptional Dehn fillings will consist of pairs $(M, \alpha)$, where $\alpha$ is exceptional.

By deleting a cusp neighborhood $N$ of $M$, we obtain a compact $3$-manifold $M\setminus N$ with torus boundary which is isotopic to $\overline{M}$. Slopes on $\partial \overline{M}$ can thus be carried by the isotopy onto $\partial(M\setminus N)\subset M$. By taking $N$ to be maximal, we define the length of a slope to be the minimum length of all curves in its isotopy class on $\partial(M\setminus N)$. 

Agol \cite{Agol:bounds} and Lackenby \cite{L:surgery} independently proved that lengths of exceptional slopes cannot exceed $6$. It follows that there are only finitely many exceptional fillings on a $1$-cusped finite-volume hyperbolic $3$-manifold. Therefore it is possible to give complete lists of exceptional Dehn fillings on $1$-cusped hyperbolic $3$-manifolds in the censuses. For censuses of at most $9$ tetrahedra, this has been done by Dunfield \cite{Dunfield:exceptional}. For the $10$-tetrahedra census, we obtained:

\begin{theorem}
  \label{th:exceptional}
  There are precisely $\numprint{439898}$ exceptional Dehn fillings on orientable $1$-cusped hyperbolic $3$-manifolds in the $10$-tetrahedra census. 
\end{theorem}

\Cref{tab:exceptional} outlines the procedure for obtaining \cref{th:exceptional}. There are $\numprint{144016}$ $1$-cusped hyperbolic $3$-manifolds in the $10$-tetrahedra census. For each of them, we enumerated all slopes with length no more than $6$, obtaining a total of $\numprint{800447}$ candidates for the exceptional census. As no duplicate slopes were enumerated for a given manifold, all candidates are distinct, and we only need to decide the hyperbolicity for each of them. 

\begin{table}[ht]
  \centering
  \begin{tabular}{crrr}
            \toprule
            Step & \multicolumn{1}{c}{Candidates} & \multicolumn{1}{c}{Discarded} & \multicolumn{1}{c}{Confirmed} \\
             \midrule
             \midrule
               Candidate generation & \numprint{800447} & 0 & 0 \\
            \midrule
              Certify hyperbolicity & \numprint{439933} & \numprint{360514} & 0\\
              \midrule
              Regina census lookup & \numprint{11540} & \numprint{24} & \numprint{428369} \\
              \midrule
              Essential $S^2$ or torus & \numprint{11} & 0 & \numprint{11529}\\
              \midrule
              Identified to be m135(1,3) & 0 & 11 & 0 \\
              \midrule
              \midrule
              Total & \numprint{800447} & \numprint{360549} & \numprint{439898}\\
              \bottomrule
    \end{tabular}
    \caption{The steps for creating the $10$-tetrahedra exceptional census}
    \label{tab:exceptional}
\end{table}

For each candidate $(M,\alpha)$, we consider its filled manifold $M(\alpha)$. Several rounds of certifying the hyperbolicity of $M(\alpha)$ with \texttt{SnapPy} reduced the remaining list by nearly a half. For Regina census lookup, we cut $M(\alpha)$ into connected sum summands and used \texttt{regina.Census.lookup} to try obtaining a name classifying it, which then tells us whether $M(\alpha)$ is hyperbolic or not. For those whose names were not obtained, we check whether they have fundamental normal surfaces which are essential $S^2$ or tori; if they do, they are not hyperbolic. Only $11$ candidates remain undetermined after this, and their filled manifolds were all identified as m135(1,3) by \texttt{SnapPy}. m135(1,3) is isometric to m412(1,-3)(-3,1) whose triangulation in \texttt{SnapPy} is certified to be geometric, hence m135(1,3) is hyperbolic.

Applying \texttt{regina.Triangulation3.isSphere} on fillings in \cref{th:exceptional} reveals that
\begin{theorem}
  \label{th:knot-exterior}
There are precisely $\numprint{1849}$ manifolds in the $10$-tetrahedra census that are exteriors of knots in $S^3$. 
\end{theorem}
For censuses of at most $9$ tetrahedra, similar work have been done by Callahan, Dean and Weeks \cite{CDW:simplest}, Champanerkar, Kofman and Patterson \cite{CKP:simplest} and Dunfield \cite{Dunfield:exceptional}.

Many conjectures and results mentioned in \cite{Dunfield:exceptional} can be tested and updated against the extended exceptional census. We plan to address these in a future publication.

\subsection{Thurston norm and twisted Alexander polynomial}

Given an orientable cusped hyperbolic $3$-manifold $M$ with Betti number $b_1(M) = 1$ and an $\SL(2,\mathbb{C})$-lift $\tilde{\rho}$ of its holonomy representation $\rho\colon \pi_1(M)\to \PSL(2,\mathbb{C})$, we consider its twisted Alexander polynomial $\tau_{2,\tilde{\rho}}(M)\in \mathbb{C}[t^{\pm 1}]$, defined by Dunfield, Friedl and Jackson \cite[\S 2]{DFJ:twisted}. It was proved in \cite{DFJ:twisted} that 
\begin{equation}
  \label{eq:twisted}
  x(M) \geq \frac{1}{2}\deg\tau_{2,\tilde\rho}(M),
\end{equation}
where $x(M)$ is the Thurston norm of a generator of $H^1(M;\mathbb{Z})$, and \cref{eq:twisted} is an equality if $M$ fibers over $S^1$. Moreover, if $M$ fibers over $S^1$, then $\tau_{2,\tilde{\rho}}(M)$ is monic.

In \cite{DGY:torsion-fibered}, Dunfield, Garoufalidis and Yoon asked whether for any orientable cusped hyperbolic $3$-manifold with $b_1(M)=1$, there is at least one lift $\tilde{\rho}$ such that \cref{eq:twisted} is an equality. An affirmative answer was then experimentally given for all such manifolds in the censuses of at most $9$ tetrahedra, a sample size of $\numprint{59068}$ manifolds.

We extended the experimentation to the $10$-tetrahedra census. There are $\numprint{143919}$ manifolds in the $10$-tetrahedra census with $b_1(M) = 1$. By applying the code in \cite{DGY:code}, we were able to compute the Thurston norms for all except the following $21$ of them: 
\small{
\begin{equation}
  \label{eq:norm-fail}
  \begin{aligned}
  o10\_077959 &&
  o10\_076513 &&
  o10\_091925 &&
  o10\_073669 &&
  o10\_067720 &&
  o10\_102779 &&
  o10\_142965 \\
  o10\_141638 &&
  o10\_146022 &&
  o10\_140310 &&
  o10\_091511 &&
  o10\_122772 &&
  o10\_147820 &&
  o10\_127911 \\
  o10\_120303 &&
  o10\_145417 &&
  o10\_136855 && 
  o10\_136895 &&
  o10\_148034 &&
  o10\_143501 &&
  o10\_138268
  \end{aligned}
\end{equation}
}

\normalsize

Out of the $21$ manifolds in \cref{eq:norm-fail}, the following $9$ were confirmed to not fiber over $S^1$, as their twisted Alexander polynomials are not monic:
\begin{equation}
  \begin{aligned}
    o10\_076513 &&
    o10\_091925 &&
    o10\_142965 &&
    o10\_091511 &&
    o10\_136855 \\
    o10\_136895 &&
    o10\_148034 &&
    o10\_143501 &&
    o10\_138268
  \end{aligned}
\end{equation}

For the remaining $\numprint{143898}$ manifolds whose Thurston norms were successfully computed, we numerically computed their twisted Alexander polynomials of all $\SL(2,\mathbb{C})$-lifts. After taking the maximum degrees, we found that \cref{eq:twisted} is an equality for all of them.

\subsection{Homology 3-spheres and the L-space conjecture}

The L-space conjecture asserts that for an irreducible $\mathbb{Q}$-homology $3$-sphere $Y$, the following are equivalent:
\begin{alphaenumerate}
  \item $\pi_1(Y)$ is left-orderable;
  \item $Y$ has non-minimal Heegaard Floer homology;
  \item $Y$ admits a taut foliation.
\end{alphaenumerate}
The equivalence of (a) and (b) was proposed by Boyer, Gordon and Watson \cite{BGW:Lspace}. The equivalence of (b) and (c) was asked by Ozsv\'ath and Szab\'o \cite{OS:holomorphic} after they proved that (c) implies (b), and upgraded to a conjecture by A. Juh\'asz \cite{Juh:HF}. For $\mathbb{Z}$-homology spheres, Calegari and Dunfield \cite{CD:lamination} proved that (c) implies (a).

Dunfield \cite{Dunfield:Lspace} showed that the L-space conjecture holds on at least $62.2\%$ of $\numprint{307301}$ hyperbolic $\mathbb{Q}$-homology $3$-spheres, obtained by Dehn fillings on orientable $1$-cusped hyperbolic $3$-manifolds in censuses of at most $9$ tetrahedra. It was conjectured in \cite{Dunfield:Lspace} that these are precisely the filled manifolds from those censuses with lengths of systoles being at least $0.2$.

While there are infinitely many hyperbolic $\mathbb{Q}$-homology $3$-spheres obtained by Dehn fillings on a finite set of orientable $1$-cusped hyperbolic $3$-manifolds, a finite list is obtained if we bound the list by the length of systoles, i.e. the shortest geodesics, from below, as the following result was shown by Hodgson and Kerckhoff \cite{HK:universal}:
\begin{lemma}[Hodgson \& Kerckhoff]
  \label{le:systole}
  Let $\overline{M}$ be a compact orientable $3$-manifold with torus boundary whose interior $M$ is hyperbolic, and $\alpha$ be a slope on $\partial\overline{M}$. If $M(\alpha)$ is hyperbolic with systole length at least $0.163$, then the normalized length of $\alpha$ is no more than $7.515$.
\end{lemma}
Here, the \textit{normalized length} of $\alpha$ is the length of $\alpha$ divided by the area of $\partial(M\setminus N)$. 

Similarly to \cref{subsec:exceptional}, for the $1$-cusped hyperbolic $3$-manifolds in the $10$-tetrahedra census, by enumerating all slopes with normalized length no more than $7.515$, excluding all exceptional slopes, and filtering by the length of systole, we obtained a list of $\numprint{1899974}$ fillings containing all pairs $(M,\alpha)$ such that $M(\alpha)$ is a hyperbolic $\mathbb{Q}$-homology $3$-sphere with length of systole at least $0.163$. Therefore,
\begin{theorem}
  \label{th:Qsphere}
  There are precisely $\numprint{1899974}$ pairs $(M,\alpha)$ where $M$ is a $1$-cusped hyperbolic $3$-manifold in the $10$-tetrahedra census and $\alpha$ is a slope such that $M(\alpha)$ is a hyperbolic $\mathbb{Q}$-homology $3$-sphere with length of systole at least $0.163$. 
\end{theorem}
Changing $\mathbb{Q}$ to $\mathbb{Z}$ in \cref{th:Qsphere}, the number of pairs becomes $\numprint{10729}$. 

Note that, although the pairs $(M,\alpha)$ in \cref{th:Qsphere} are all distinct, the filled manifolds $M(\alpha)$ are not; for instance, there are at least $22$ pairs in the list such that $M(\alpha)$ is identical to m007(3,1). It remains to group the pairs into distinct isometry classes by $M(\alpha)$, compare the grouped list against the list of $\numprint{307301}$ manifolds in \cite{Dunfield:Lspace}, and test the L-space conjecture on it. We wish to address these in a future publication. 

\subsection{Closed totally geodesic surfaces}

An embedded surface $S$ in a hyperbolic $3$-manifold $M$ is called \textit{totally geodesic} if every geodesic in $S$ is also a geodesic in $M$. The existence of closed totally geodesic surfaces is of particular interest.

An algorithm that is guaranteed to find all closed totally geodesic surfaces (up to isotopy) in a given orientable cusped finite-volume hyperbolic $3$-manifold was built and used to find several examples of cusped hyperbolic $3$-manifolds with closed totally geodesic surfaces by Basilio, Lee and Malionek \cite{BLM:geodesic}; however, none was found in the censuses of at most $9$ tetrahedra. By applying their algorithm to the $10$-tetrahedra census, we found that 
\begin{theorem}
  The $2$-cusped manifold o10\_143602, whose volume is approximately $9.1345$ and first homology group is $\mathbb{Z}^2\oplus \mathbb{Z}_6$, 
  is the only manifold containing closed totally geodesic surfaces in the orientable cusped hyperbolic censuses of at most $10$ tetrahedra. Moreover, o10\_143602 contains exactly one such surface up to isotopy, which is nonorientable with Euler characteristic $-1$. 
\end{theorem}
In particular, among all orientable cusped hyperbolic $3$-manifolds containing closed totally geodesic surfaces, o10\_143602 is the one whose minimal triangulations have the fewest tetrahedra; in other words, it is the simplest such example.


\bibliography{references}

\appendix
\section{Miscellaneous tables}
\label{appen:tables}
\subsection{Summary of tests}
\label{subappen:summary}
\begin{table}[ht!]
\centering
\begin{tabular}{ccc}
          \toprule
          Certification & Name  & Function(s)\\
          \midrule
          \multirow{2}{*}{Nonminimality} & Greedy nonminimality & \texttt{regina.Triangulation3.intelligentSimplify}\\
          \cmidrule{2-3}
          & Exhaustive nonminimality & \texttt{regina.Triangulation3.simplifyExhaustive}\\
          \midrule
          \multirow{3}{*}{\Centerstack{Nonminimality \\ or \\
          nonhyperbolicity}} & \multirow{3}{*}{Special surfaces} & \texttt{regina.NormalSurfaces}\\
          && \texttt{regina.NormalSurface.isTwoSided}\\
          && \texttt{regina.NormalSurface.isSolidTorus}\\
          \midrule
          Hyperbolicity & Certify hyperbolicity & \texttt{snappy.Manifold.verify\_hyperbolicity}\\
          \bottomrule
\end{tabular}
\caption{Summary of tests, where functions are given in python convention}
\label{tab:tests}
\end{table}

\subsection{Manifold counts by cusp number}
\label{subappen:cusps}

\begin{table}[H]
  \centering
  \begin{tabular}{crrrrrrrrr}
            \toprule
            & \multicolumn{9}{c}{Tetrahedra} \\
            \midrule
            Cusps & 2 & 3 & 4 & 5 & 6 & 7 & 8 & 9 & 10 \\
            \midrule
            \midrule
            1 & 2 & 9 & 52 & 223 & 913 & \numprint{3388} & \numprint{12241} & \numprint{42279} & \numprint{144016}\\
            2 & 0 & 0 & 4 & 11 & 48 & 162 & \numprint{591} & \numprint{1934} & \numprint{6585}\\
            3 & 0 & 0 & 0 & 0 & 1 & 2 & 13 & \numprint{36} & \numprint{123}\\
            4 & 0 & 0 & 0 & 0 & 0 & 0 & 1 & \numprint{1} & \numprint{5}\\
            5 & 0 & 0 & 0 & 0 & 0 & 0 & 0 & \numprint{0} & \numprint{1}\\
            \bottomrule
  \end{tabular}
  \caption{Orientable cusped hyperbolic $3$-manifold counts by cusp and tetrahedra number}
  \label{tab:cusps}
  \end{table}

In \cite{AS:minimum}, Adams and Sherman proved that the least number of tetrahedra required to build a $5$-cusped hyperbolic $3$-manifold is $10$, with an example being the complement of the link $L10a174$. It follows from \cref{tab:cusps} that the complement of $L10a174$, named as o10\_150729 in the $10$-tetrahedra census, is the unique orientable $5$-cusped hyperbolic $3$-manifold admitting an ideal triangulation of $10$ tetrahedra. Furthermore, o10\_150729 has exactly $2$ combinatorially distinct minimal ideal triangulations.

\end{document}